\documentclass[12pt]{article}
\usepackage{amsfonts}
\usepackage{theorem}

\usepackage[OT2,OT1]{fontenc}
\newcommand\cyr{%
\renewcommand\rmdefault{wncyr}%
\renewcommand\sfdefault{wncyss}%
\renewcommand\encodingdefault{OT2}%
\normalfont
\selectfont}
\DeclareTextFontCommand{\textcyr}{\cyr}

\newtheorem{theorem}{Theorem}
\newtheorem{definition}{Definition}
\newtheorem{corollary}{Corollary}
\newtheorem{lemma}{Lemma}
\newtheorem{remark}{Remark}
\newtheorem{example}{Example}
\usepackage[russian, english]{babel}
\usepackage{centernot}
\usepackage[unicode,colorlinks,plainpages=false]{hyperref}
\newcommand{\openbox}{$\begin{array}{c}
\hspace*{-0.55em}\sqcap \hspace*{-0.60em}\\[-0.4em] \hline
\multicolumn{1}{c}{\hspace*{-0.60em}}\\[-0.8em]
\end{array}$}

\begin{document}

\centerline{\bf Separators of Ideals in Multiplicative Semigroups}
\centerline{\bf of Unique Factorization Domains}
\medskip

\centerline{Attila Nagy}

\bigskip

\small
\begin{abstract}
In this paper we show that if $I$ is an ideal of a commutative semigroup $C$ such that the separator $SepI$ of $I$ is not empty then the factor semigroup $S=C/P_I$ ($P_I$ is the principal congruence on $C$ defined by $I$) satisfies  Condition $(*)$: $S$ is a commutative monoid with a zero; The annihilator $A(s)$ of every non identity element $s$ of $S$ contains a non zero element of $S$; $A(s)=A(t)$ implies $s=t$ for every $s, t\in S$. Conversely, if $\alpha$ is a congruence on a commutative semigroup $C$ such that the factor semigroup $S=C/\alpha$ satisfies Condition $(*)$ then there is an ideal $I$ of $C$ such that $\alpha =P_I$. Using this result for the multiplicative semigroup $D_{mult}$ of a unique factorization domain $D$, we show that $P_{J(m)}=\tau _m$ for every nonzero element $m\in D$, where $J(m)$ denotes the ideal of $D$ generated by $m$, and $\tau _m$ is the relation on $D$ defined by $(a, b)\in \tau _m$ if and only if $gcd(a, m)\sim gcd(b, m)$ ($\sim$ is the associate congruence on $D_{mult}$). We also show that if $a$ is a nonzero element of a unique factorization domain $D$ then $d(a)=|D'/P_{J([a])}|$, where $d(a)$ denotes the number of all non associated divisors of $a$, $D'=D/\sim$, and $[a]$ denotes the $\sim$-class of $D_{mult}$ containing $a$.
As an other application, we show that if $d$ is one of the integers $-1$, $-2$, $-3$, $-7$, $-11$, $-19$, $-43$, $-67$, $-163$ then, for every nonzero ideal $I$ of the ring $R$ of all algebraic integers of an imaginary quadratic number field ${\mathbb Q}[\sqrt d]$, there is a nonzero element $m$ of
$R$ such that $P_I=\tau _m$.


\end{abstract}

\normalsize

\medskip
\section{Introduction}

Let $S$ be a commutative semigroup with a zero element $0$. If $s$ is an arbitrary element of $S$, then the annihilator $A(s)=\{ x\in S:\ xs=0\}$ is an ideal of $S$. An element $s\in S$ is said to be a torsion element if $A(s)\neq \{ 0\}$. Let $S_T$ denote the set of all torsion elements of $S$.

\begin{lemma}\label{lm1} If $S$ is a commutative semigroup with a zero such that $|S|\geq 2$, then $S_T$ is an ideal of $S$.
\end{lemma}

\noindent
{\bf Proof}. Let $S$ be a commutative semigroup with a zero element $0$ such that $|S|\geq 2$. Then $0\in S_T$. If $a\in S_T$, then there is an element $b\neq 0$ of $S$ such that $b\in A(a)$. Thus, for an arbitrary element $s\in S$, we have $(sa)b=s(ab)=s0=0$ which implies
$0\neq b\in A(sa)$. Hence $sa\in S_T$.
\hfill\openbox

\begin{definition}\label{df1} We shall say that a semigroup $S$ satisfies Condition $(*)$ if the following three conditions hold.
\begin{enumerate}
\item $S$ is a commutative monoid with a zero;
 \item Every non identity element of $S$ is a torsion element;
 \item $A(s)=A(t)$ implies $s=t$ for every $s, t\in S$.
\end{enumerate}
\end{definition}

\medskip

We note that a one-element semigroup $S=\{ e\}$ satisfies Condition $(*)$ of Definition~\ref{df1} because $e$ is both an identity element and a zero element of $S$, and $S$ has no non identity element.

\begin{example}\label{ex1}
Let $S=\{ 0, 1, 2\}$ be the semigroup defined by Table~\ref{Fig1}.
\begin{table}[ht]
\center{\begin{tabular}{l|l l l }
 &$1$&$2$&$0$\\ \hline
$1$&$1$&$2$&$0$\\
$2$&$2$&$0$&$0$\\
$0$&$0$&$0$&$0$\\
 \end{tabular}}
\caption{}\label{Fig1}
\end{table}

This semigroup is a commutative monoid with a zero in which $1$ is the identity element and $0$ is the zero element. Since $A(1)=\{ 0\}$, $A(2)=\{ 2, 0\}$,
$A(0)=\{ 1, 2, 0\}$ then $S$ satisfies Condition $(*)$ of Definition~\ref{df1}.
\end{example}

\medskip

The natural quasi-order on a monoid S is the relation defined by divisibility: $a\leq b$ for some $a, b\in S$ if $bS\subseteq aS$.
If this quasi-order is a partial order on $S$, then we speak of natural partial order and call the monoid $S$ a naturally partially
ordered monoid. In such a monoid, the identity element 1 is the unique invertible element which means that, for arbitrary $a, b\in S$, $ab=1$ if and only if $a=1$ and $b=1$.

\begin{lemma}\label{lm2} Every semigroup satisfying Condition $(*)$ of Definition~\ref{df1} is a naturally partially ordered monoid.
\end{lemma}

\noindent
{\bf Proof}. Let $S$ be a semigroup satisfying condition $(*)$ of Definition~\ref{df1}. Assume $a\leq b$ and $b\leq a$ for some $a, b\in S$. Then
$Sa=Sb$ and so there are elements $x, y\in S$ such that $a=xb$ and $b=ya$. Let $s\in A(a)$ be arbitrary. Then $bs=(ya)s=y(as)=y0=0$ and so $s\in A(b)$.
Hence $A(a)\subseteq A(b)$. Similarly, $A(b)\subseteq A(a)$ and so $A(a)=A(b)$. As $S$ satisties $(*)$, we get $a=b$. Thus the preorder $\leq$ is also antisymmetric. Consequently, $\leq$ is a partial order and so $S$ is a naturally partially ordered monoid.\hfill\openbox

\section{Commutative semigroups}

Let $S$ be a semigroup. By the idealizer $IdA$ of a subset $A$ of $S$ we mean the set of all elements $x$ of $S$ with the property $xA\subseteq A$, $Ax\subseteq A$ (\cite{kantorovics}). The idealizer of the empty set is $S$. It is clear that, for every subset $A$ of a semigroup $S$, the idealizer of $A$ is either empty or a subsemigroup of $S$. The intersection $SepA=IdA\cap Id(S\setminus A)$ is called the separator of $A$ (\cite{nagysep}). By the above remark, the separator of a subset of a semigroup $S$ is either empty or a subsemigroup of $S$. It is clear that $SepA=Sep(S\setminus A)$. Especially, $Sep\emptyset =SepS=S$.
If $A$ is a subsemigroup of a semigroup $S$, then $A\cup SepA$ is also a subsemigroup of $S$ by Theorem 2 of \cite{nagysep}. In Theorem 3 of \cite{nagysep} it is proved that $SepA\subseteq A$ or $SepA\subseteq S\setminus A$ for every subset $A$ of a semigroup $S$.
An ideal $I$ of a semigroup $S$ is called a proper ideal if $I\neq S$.
By Remark 3 of \cite{nagysep}, if $I$ is a proper ideal of a semigroup $S$, then $SepI\subseteq S\setminus I$.

\medskip

Using the notation of \cite{clifford2}, if $H$ is a subset of a semigroup $S$, then let
\[H\dots a=\{(x,y)\in S\times S:\ xay\in H\}, \quad a\in S\]
and
\[P_H=\{(a,b)\in S\times S:\ H\dots a=H\dots b\}.\]
It is easy to see that $P_H$ is a congruence on $S$; this congruence is called the principal congruence on $S$ defined by $H$ (\cite{clifford2}).

For notations and notions not defined here, we refer to \cite{clifford1}, \cite{clifford2} and \cite{nagybook}.

\medskip

\begin{theorem}\label{th1} Let $C$ be a commutative semigroup. If $I$ is an ideal of $C$ and $SepI\neq \emptyset$, then $P_I$ is a congruence on $C$ such that $I$ and $SepI$ are $P_I$-classes of $C$, and the factor semigroup $S=C/P_I$ satisfies Condition $(*)$ of Definition~\ref{df1}.

Conversely, if $\alpha$ is a congruence on a commutative semigroup $C$ such that the factor semigroup $S=C/\alpha$ satisfies Condition $(*)$ of Definition~\ref{df1} then there is an ideal $I$ of $C$ such that $\alpha =P_I$.
\end{theorem}

\noindent
{\bf Proof}. Let $C$ be an arbitrary commutative semigroup and $I$ an arbitrary ideal of $C$ with condition $SepI\neq \emptyset$.

If $I=C$ then $SepI=C$ and $P_I$ is the universal relation on $C$; in this case $I$ and $SepI$ are $P_I$-classes of $C$, and the factor semigroup $S=C/P_I$ is a one-element semigroup and so it satisfies Condition $(*)$ of Definition~\ref{df1}.

Next, consider the case when $I\neq C$. By Remark 3 of \cite{nagysep},
$I\cap SepI=\emptyset$. By Theorem 2 and Theorem 3 of \cite{nagycommon}, $P_I$ is a congruence on $C$ such that $SepI$ is the identity element of the factor semigroup $C/P_I$, and $I$ is a union of $P_I$-classes of $C$. As $I...a=C\times C$ for every $a\in I$, the ideal $I$ is a $P_I$-class of $C$; $I$ is the zero of $C/P_I$. Thus $C/P_I$ is a commutative monoid with a zero element.

For an arbitrary element $x$ of $C$, let $[x]_I$ denote the $P_I$-class of $C$ containing the element $x\in C$.
If $[c]_I$ is a non identity element of $C/P_I$, then $c\in C\setminus SepI$. As $c\in C=IdI$, we have $c\notin Id(C\setminus I)$. Thus there is an element $b\in (C\setminus I)$ such that $cb\in I$. Then $[b]_I$ is a nonzero element of $C/P_I$ such that $[c]_I[b]_I$ is the zero of $C/P_I$. Hence $[b]_I\in A([c]_I)$. Consequently every non identity element of $C/P_I$ is a torsion element.

Assume $[a]_I\neq [b]_I$ for some $a, b\in C$. It means that $(a, b)\notin P_I$ and so \[I...a\neq I...b.\] Thus there is a couple $(x, y)\in C\times C$ such that $(x, y)\in I...a$ and $(x, y)\notin I...b$ or $(x, y)\in I...b$ and $(x, y)\notin I...a$. Consider the case \[(x, y)\in I...a\quad \hbox{and}\quad (x, y)\notin I...b.\] Then \[xay\in I\quad xby\notin I\] from which it follows that \[axy\in I,\quad \hbox{and}\quad bxy\notin I.\] Thus \[[a]_I[xy]_I=0\quad \hbox{and}\quad [b]_I[xy]_I\neq 0\] in the factor semigroup $S/P_I$. Hence  \[[xy]_I\in A([a]_I),\quad [xy]_I\notin A([b]_I).\] Consequently $A([a]_I)\neq A([b]_I)$.

We also get $A([a]_I)\neq A([b]_I)$ in case $(x, y)\in I...b$ and $(x, y)\notin I...a$. Thus, in both cases, $A([a]_I)\neq A([b]_I)$. Then the assumption $A([a]_I)=A([b]_I)$ implies $[a]_I=[b]_I$ for every $a, b\in C$.

Summarizing our above results, we can see that the factor semigroup $S=C/P_I$ has all three properties contained Condition $(*)$ of Definition~\ref{df1}.

To prove the converse, let $\alpha$ be a congruence on the semigroup $C$ such that the factor semigroup $S=C/\alpha$ satisfies Condition $(*)$ of Definition~\ref{df1}. Let $I$ and $H$ denote the $\alpha$-classes of $C$ which are the zero element and the identity element of $C/\alpha$, respectively. It is clear that $I$ is an ideal of $C$.

We show that $SepI=H$. Let $[x]_{\alpha}$ denote the $\alpha$-class of $C$ containing $x\in C$. As $H$ is the identity element of $C/\alpha$, $H[x]_{\alpha}\subseteq [x]_{\alpha}$ for every $x\in C$. Then $H\subseteq SepI$. Let $a\notin H$ be an arbitrary element of $C$. Then $[a]_{\alpha}$ is not the identity element of $C/\alpha$ and so $A([a]_{\alpha})\neq \{0\}$ in $C/\alpha$, because the factor semigroup $S=C/\alpha$ satisfies condition $(*)$ of Definition~\ref{df1}. Thus there is an element $b\notin I$ of $C$ such that $[a]_{\alpha}[b]_{\alpha}\subseteq I$ and so $ab\in I$. Consequently
$a\notin Id(C\setminus I)$ and so $a\notin SepI$. From this result it follows that $H=SepI$.

We show that $\alpha =P_I$. Assume $(a, b)\in \alpha$ for some $a, b\in C$. Then, for every $x, y\in C$, there is an element $c\in C$ such that
\[xay, xby\in [c]_{\alpha}.\] As $I$ is an $\alpha$-class of $C$, $xay\in I$ if and only if $xby\in I$. Hence \[I...a=I...b,\] that is, $(a,b)\in P_I$. Thus \[\alpha \subseteq P_I.\]

To show $P_I\subseteq \alpha$, assume $(a,b)\in P_I$ for some
$a, b\in C$. If $(a,b)\notin \alpha$ then $[a]_{\alpha}\neq [b]_{\alpha}$ and so $A([a]_{\alpha})\neq A([b]_{\alpha})$, because the factor semigroup $S=C/\alpha$ satisfies Condition $(*)$ of Definition~\ref{df1}. From this it follows that there is an element $[x]_{\alpha}\in C/\alpha$ such that $[a]_{\alpha}[x]_{\alpha}\subseteq I$ and $[b]_{\alpha}[x]_{\alpha}\subseteq (C\setminus I)$ or
$[a]_{\alpha}[x]_{\alpha}\subseteq (C\setminus I)$ and $[b]_{\alpha}[x]_{\alpha}\subseteq I$, because $I$ is the zero of $C/\alpha$. Consider the case
\[[a]_{\alpha}[x]_{\alpha}\subseteq I\quad \hbox{and}\quad [b]_{\alpha}[x]_{\alpha}\subseteq (C\setminus I).\]
Then \[ax\in I\quad \hbox{and}\quad bx\notin I\] and so, for an arbitrary
$e\in H=SepI$, \[eax\in I\quad \hbox{and}\quad ebx\notin I\] which imply \[(e, x)\in I...a\quad \hbox{and}\quad (e, x)\notin I...b.\] Consequently \[I...a\neq I...b.\] We get the same result in the other case.at case when $[a]_{\alpha}[x]_{\alpha}\subseteq (C\setminus I)$ and $[b]_{\alpha}[x]_{\alpha}\subseteq I$. Hence \[(a, b)\notin P_I\] which is a contradiction. Consequently $P_I\subseteq \alpha$ and so $\alpha =P_I$.\hfill\openbox

\begin{remark}\label{rmttt} Let $C$ be a commutative semigroup with a zero $0$ such that the nonzero elements of $C$ form a subsemigroup $C^*$ of $C$. For an ideal $I\neq \{ 0\}$ of $C$, let $I^*=I\setminus \{ 0\}$. Then $I^*$ is an ideal of $C^*$ and the restriction of $P_I$ to $C^*$ equals $P_{I^*}$, where $P_I$ denotes the main congruence on $C$ defined by $I$ and  $P_{I^*}$ is the main congruence on $C^*$ defined by $I^*$. As $I$ is a $P_I$-class and $I^*$ is a $P_{I^*}$-class by Theorem~\ref{th1}, we have $C/P_I\cong C^*/P_{I^*}$.
\end{remark}

\bigskip

A proper ideal $I$ of a semigroup $S$ is called a prime ideal if, for every $a, b\in S$, the assumption $ab\in I$ implies $a\in I$ or $b\in I$.
It is known that a proper ideal $I$ of a commutative semigroup $S$ is prime if and only if, for every ideal $A$ and $B$ of $S$, the assumption $AB\subseteq I$ implies $A\subseteq I$ or $B\subseteq I$.

\medskip

A proper ideal $I$ of a semigroup $S$ is called a maximal ideal of $S$ if, for every ideal $A$ of $S$, the assumption $I\subseteq A\subseteq S$ implies $I=A$ or $A=S$.



\begin{theorem}\label{th2} On a maximal ideal $M$ of a commutative semigroup $S$, the following assertions are equivalent.
\begin{enumerate}
\item[(i)] $M$ is a prime ideal.
\item[(ii)] $SepM\neq \emptyset$.
\item[(iii)] The factor semigroup $S/P_M$ is a two-element monoid with a zero.
\end{enumerate}
\end{theorem}

\noindent
{\bf Proof}. $(i)$ implies $(ii)$: Let $M$ be a prime ideal. Then $SepM=S\setminus M\neq \emptyset$ by Theorem 9 of \cite{nagysep}.

$(ii)$ implies $(i)$: Let $M$ be a maximal ideal of $S$ such that $SepM\neq \emptyset$. Assume $AB\subseteq M$  for some ideals $A\not \subseteq M$,
$B\not \subseteq M$ of $S$. As $M$ is a maximal ideal of $S$ and both of $M\cup A$, $M\cup B$ are ideals of $S$, we have $M\cup A=S$, $M\cup B=S$. Thus $(S\setminus M)\subseteq A$ and $(S\setminus M)\subseteq B$. As $SepM\subseteq (S\setminus M)$, we have $SepM\subseteq A$, $Sep M\subseteq B$. Then, for arbitrary elements $x, y\in SepM$, $xy\in AB\subseteq M$. This is a contradiction, because $SepM$ is a subsemigroup of $S$ and so $xy\in SepM\subseteq (S\setminus M)$.
Consequently, for every ideals $A$ and $B$ of $S$, the assumption $AB\subseteq M$ implies $A\subseteq M$ or $B\subseteq M$. Hence $M$ is a prime ideal, because $S$ is a commutative semigroup.

$(i)$ implies $(iii)$: Let $M$ be a prime ideal of $S$. Then $S\setminus M$ is a subsemigroup of $S$. By Remark 1, Theorem 9 and Theorem 8 of \cite{nagysep}, $SepM=Sep(S\setminus M)=(S\setminus M)$. Thus the congruence $P_M$ on $S$ has two classes by Theorem~\ref{th1}. Hence the factor semigroup is a two-element monoid with a zero.

$(iii)$ implies $(i)$: Assume that the factor semigroup $S/P_M$ is a two-element monoid with a zero. Let $c\in (S\setminus M)$ be an element
with $M...c=S\times S$. As $K=\{ s\in S: M...s=S\times S\}$ is an ideal of $S$ such that $M\subset K$ (because $c\in K$ and $c\notin M$), we get $K=S$
by the condition that $M$ is a maximal ideal of $S$. Hence $P_M$ is the universal relation on $S$. This is a contradiction. From this it follows that $M...c\neq S\times S$ for arbitrary $c\in (S\setminus M)$. As $M...s=S\times S$ for every $s\in M$, we have that $M$ is a $P_M$-class of $S$. As the factor semigroup is a two-element monoid with a zero, $(S\setminus M)$ is a subsemigroup of $S$ which is the identity element of the factor semigroup $S/P_M$. Then $(S\setminus M)=Sep(S\setminus M)$ and so $M$ is a prime ideal by Theorem 8 and Theorem 9 of \cite{nagysep}.\hfill\openbox

\medskip



\section{Multiplicative semigroups of unique factorization domains}

Let $D$ be an integral domain with an identity element $e$. Then the multiplicative semigroup $D_{mult}$ of $D$ is a commutative monoid and so every ideal of $D_{mult}$ has a non empty separator. For arbitrary $a, b\in D$, let $a\sim b$ denotes the fact that $a$ and $b$ are associated elements. It is known that $\sim$ is an equivalence relation on $D$. It is easy to see that $\sim$ is a congruence on $D_{mult}$, and $D'=D_{mult}/\sim$ is a commutative monoid with the identity element $[e]_{\sim}$ ($[x]_{\sim}$ denotes the $\sim$-class of $D$ containing $x\in D$). $[e]_{\sim}$ is the only unit of $D'$, and $[0]_{\sim}$ is the only divisor of zero.

\begin{lemma}\label{simsub} Let $D$ be an integral domain with an identity element. Then, for an arbitrary ideal $I$ of $D_{mult}$, we have $\sim \subseteq P_I$.
\end{lemma}

\noindent
{\bf Proof}. Let $a$ and $b$ be arbitrary elements of $D$ with $a\sim b$. Then there are units $\epsilon$, $\epsilon '$ in $D$ such that $a=\epsilon b$ and $b=\epsilon 'a$. Assume $xay\in I$ for some $x, y\in D$. Then $xby=\epsilon 'xay\in I$. Similarly, $xby\in I$ implies $xay=\epsilon xby\in I$. Thus $I...a=I...b$ and so $(a,b)\in P_I$. Hence $\sim \subseteq P_I$. \hfill\openbox

\begin{lemma}\label{phisep} Let $D$ be an integral domain with an identity element and $\varphi$ the canonical homomorphism of $D_{mult}$ onto $D'=D_{mult}/\sim$. For an arbitrary ideal $I$ of $D_{mult}$, $\varphi (SepI)=Sep(\varphi (I))$.
\end{lemma}

\noindent
{\bf Proof}. Let $I$ be an arbitrary ideal of $D_{mult}$. By Theorem~\ref{th1}, $I$ and $SepI$ are $P_I$-classes of $D_{mult}$. Using Lemma~\ref{simsub}, we can see that $\varphi (SepI)=Sep(\varphi (I))$.\hfill\openbox

\bigskip

In the next, $D$ will denote a unique factorization domain.
For arbitrary elements $a$ and $b$ of $D$, there exist the greatest common divisor $gcd(a, b)$ of $a$ and $b$ in $D$. It is clear that, in $D_{mult}/\sim$, we have $[gcd(a, b)]_{\sim}=gcd([a]_{\sim}, [b]_{\sim})$ for arbitrary $a, b\in D$.

For arbitrary $m\in D$, let $\tau _m$ denote the relation on $D$ as follows: for elements $a$ and $b$ of $D$,
$(a, b)\in \tau _m$ if and only if the greatest common divisors $gcd(a, m)$ and $gcd(b, m)$ are associated, that is, $gcd(a, m)\sim gcd(b, m)$ . As the relation $\sim$ is an equivalence relation on $D$, $\tau _m$ is also an equivalence relation on $D$.

\begin{lemma}\label{two} In an arbitrary unique factorization domain $D$, $\tau _0=\sim$.
\end{lemma}

\noindent
{\bf Proof}. As $gcd(a, 0)\sim a$ for an arbitrary $a\in D$, we have that, for every $a, b\in D$, $(a, b)\in \tau _0$ if and only if $a\sim b$. Thus $\tau _0=\sim$.\hfill\openbox

\medskip

For an arbitrary element $m$ of a unique factorization domain $D$, let $J(m)$ denote the ideal of $D$ generated by $m$. It is known that $J(m)=mD$. This ideal plays an important role in the description of the relation $\tau _m$ for an arbitrary non zero element $m$ of $D$.

\begin{theorem}\label{th3} Let $m$ be an arbitrary nonzero element of a unique factorization domain $D$. Then $\tau _m=P_{J(m)}$ and so $\tau _m$ is a congruence on $D_{mult}$ such that the factor semigroup $D_{mult}/\tau _m$ satisfies Condition $(*)$ of Definition~\ref{df1}.
\end{theorem}

\noindent
{\bf Proof}. Let $D$ be a unique factorization domain and $m\in D\setminus \{ 0\}$ be arbitrary.

First consider the case when $m \sim e$, that is, $m$ is a unit of $D$. Then $J(m)=mD=D$ and so $Sep(J(m))=D$. Thus both of $P_{J(m)}$ and $\tau _m$ are equal to the universal relation on $D$; in this case our assertions hold.

Next, consider the case when $m\not\sim e$. Then $J(m)\neq D$. As $D$ has an identity element, $Sep(J(m))\neq \emptyset$.
By Theorem~\ref{th1}, the factor semigroup $D_{mult}/P_{J(m)}$ satisfies Condition $(*)$ of Definition~\ref{df1}. $J(m)$ and $Sep(J(m))$ are different $P_{J(m)}$-classes of $D_{mult}$; $J(m)$ is the zero element, $Sep(J(m))$ is the identity element of the factor semigroup $D_{mult}/P_{J(m)}$.

We prove $\tau _m=P_{J(m)}$ by showing that $\tau _m\subseteq P_{J(m)}$ and, for every $a, b\in D$, $(a, b)\notin \tau _m$ implies $(a, b)\notin P_{J(m)}$.

Let $a$ and $b$ be arbitrary elements of $D$. Let \[d_a=gcd(a, m)\quad \hbox{and}\quad d_b=gcd(b, m).\]
Then there are elements $a^*, b^*, m', m''\in D$ such that
\[a=d_aa^*,\quad m=d_am',\quad b=d_bb^*,\quad m=d_bm''.\]
We note that \[gcd(a^*, m')\sim e\quad \hbox{and}\quad gcd(b^*, m'')\sim e.\]

Assume \[(a, b)\in \tau _,\] that is, $d_a \sim d_b$. We show that $(a, b)\in P_{J(m)}$.
As $d_a | d_b$ and $d_b | d_a$, there are elements $x, y\in D$ such that \[d_b=d_ax,\quad d_a=d_by.\]  Assume
\[uav\in J(m)\] for some $u, v\in D$. Then \[uav=mt\] for some $t\in D$ and so \[ud_aa^*v=d_am't.\] As $d_a\neq 0$, we get \[ua^*v=m't\] which implies
\[m' | uva^*.\] As $gcd(a^*, m')\sim e$, we have \[m' | uv.\] Then \[uv=m'k\] for some $k\in D$ and so
\[ubv=buv=bm'k=d_bb^*m'k=d_axb^*m'k=mxb^*k\in J(m).\]
We can prove in a similar way that $ubv\in J(m)$ implies $uav\in J(m)$. Consequently \[J(m)...a=J(m)...b\] and so \[(a, b)\in P_{J(m)}.\] Hence \[\tau _m\subseteq P_{J(m)}.\]

Next, consider the case when \[(a, b)\notin \tau _m,\] that is, $d_a\not\sim d_b$. Then $d_a \centernot| d_b$ or $d_b \centernot| d_a$. We show that
$(a, b)\notin P_{J(m)}$.

Assume \[d_a \centernot| d_b.\] Then $m'' \centernot| m'$, because $m''|m'$ would imply $m'=m''w$ ($w\in D$) from which we would get $d_am''w=d_am'=m=d_bm''$. As $m''\neq 0$, this last equation implies $d_aw=d_b$ and so we would get $d_a|d_b$ which contradicts our assumption $d_a \centernot| d_b$.

From \[m'ae=m'a =m'a^*d_a=ma^*\in J(m)\] it follows that \[(m', e)\in J(m)...a.\] We show that $(m', e)\notin J(m)...b$. Assume, in an indirect way, that $(m', e)\in J(m)...b$. Then $m'b\in J(m)$ and so $m'b=mw$ for some $w\in D$. Thus \[m'b^*d_b=m''d_bw\] and so \[m'b^*=m''w,\] because $d_b\neq 0$. Hence \[m''|m'b^*.\] As
$gcd(m'', b^*)\sim e$, we have $m''|m'$ which is a contradiction. Consequently \[(m', e)\notin J(m)...b.\] Hence \[J(m)...a\neq J(m)...b.\]

In case $d_b \centernot| d_a$, we can prove (as above) that $(m'', e)\in J(m)...b$ and $(m'', e)\notin J(m)...a$. Thus $J(m)...a\neq J(m)...b$.
In both cases, \[(a, b)\notin P_{J(m)}.\] This result together with $\tau _m\subseteq P_{J(m)}$ imply \[\tau _m=P_{J(m)}.\]
Using Theorem~\ref{th1}, we get that $\tau _m$ is a congruence on $D_{mult}$ such that the factor semigroup
$D_{mult}/\tau _m$ satisfies Condition $(*)$ of Definition~\ref{df1}.\hfill\openbox

\medskip

\begin{corollary}\label{ppp} If $m$ is a nonzero and non unit element of a unique factorization domain $D$ then $Sep(J(m))=\{ k\in D:\ gcd(k, m)\sim e\}$.
\end{corollary}

\noindent
{\bf Proof}. By Theorem~\ref{th3}, $P_{J(m)}=\tau _m$. By Theorem~\ref{th1}, $Sep(J(m))$ is a $P_{J(m)}$-class, and so a $\tau _m$-class of $D$.
As $e\in Sep(J(m))$,  we have $Sep(J(m))=\{ k\in D:\ gcd(k, m)\sim gcd(e, m)\sim e\}$.\hfill\openbox

\begin{corollary}\label{cr1} Let $m$ be an arbitrary element of a unique factorization domain $D$. For arbitrary elements $a, b, s\in D$, $gcd(a, m)\sim gcd(b, m)$ implies $gcd(sa, m)\sim gcd(sb, m)$.
\end{corollary}

\noindent
{\bf Proof}. By Lemma~\ref{two} and Theorem~\ref{th3}, $\tau _m$ is a congruence on $D_{mult}$ for every element $m\in D$. Thus our assertion is obvious.\hfill\openbox

\medskip

By Lemma~\ref{simsub}, for an arbitrary ideal $I$ of the multiplicative semigroup $D_{mult}$ of an integral domain $D$ with an identity element, we have $\sim \subseteq P_I$. Then the relation $P_I/\sim$ on the factor semigroup $D'=D_{mult}/\sim$ defined by
$([a]_{\sim}, [b]_{\sim})\in P_I/\sim$ for some $a, b\in D_{mult}$ if and only if $(a, b)\in P_I$ is a congruence on $D_{mult}/\sim$ and, by Theorem 5.6 of \cite{Howie}, $(D_{mult}/\sim)/(P_I/\sim)\cong D_{mult}/P_I$.

\begin{lemma}\label{ttt} Let $D$ be an integral domain with an identity element. In the factor semigroup $D'=D_{mult}/\sim$, we have $P_{J(m)}/\sim=P_{J([m]_{\sim})}$ for every nonzero element $m$ of $D$, where $J([m]_{\sim})$ denotes the ideal of $D'$ generated by $[m]_{\sim}=\varphi (m)$, and $P_{J([m]_{\sim})}$ is the principal congruence on $D'$ defined by $J([m]_{\sim})$.
\end{lemma}

\noindent
{\bf Proof}. As $J(m)$ is a union of $\sim$-classes of $D_{mult}$ by Lemma~\ref{simsub}, $[x]_{\sim}[a]_{\sim}[y]_{\sim}\in \varphi (J(m))$ if and only if $xay\in J(m)$. Thus $([a]_{\sim}, [b]_{\sim})\in P_{\varphi (J(m))}$ if and only if $(a,b)\in P_{J(m)}$, that is, $([a]_{\sim}, [b]_{\sim})\in P_{J(m)}/\sim$. Consequently $P_{J(m)}/\sim=P_{\varphi (J(m))}$. As $\varphi (J(m))=J([m]_{\sim})$, the lemma is proved.\hfill\openbox

\medskip

If $D$ is a unique factorization domain then, for arbitrary elements $[a]_{\sim}$ and $[b]_{\sim}$ of the factor semigroup $D'=D_{mult}/\sim$, $[a]_{\sim}$ and $[b]_{\sim}$ are associated if and only if $[a]_{\sim}=[b]_{\sim}$. For an element $[m]_{\sim}\in D'$, let $\tau '_{[m]_{\sim}}$ denote the relation on $D'$ defined by: $([a]_{\sim}, [b]_{\sim})\in \tau '_{[m]_{\sim}}$ for some elements $[a]_{\sim}$ and $[b]_{\sim}$ of $D'$ if and only if $gcd([a]_{\sim}, [m]_{\sim})=gcd([b]_{\sim}, [m]_{\sim})$. As $\tau _0=\sim$ by Lemma~\ref{two}, and $\sim \subseteq P_{J(m)}=\tau _m$ for every $0\neq m\in D$ by Theorem~\ref{th3} and Lemma~\ref{simsub}, the congruence $\tau _m/\sim$ exists.

\begin{lemma}\label{sss} Let $D$ be a unique factorization domain. In the factor semigroup $D'=D_{mult}/\sim$, we have $\tau _m/\sim=\tau '_{[m]_{\sim}}$.
\end{lemma}

\noindent
{\bf Proof}. For arbitrary $a, b\in D_{mult}$, $([a]_{\sim}, [b]_{\sim})\in \tau _m/\sim$, that is, $(a, b)\in \tau _m$ iff
$gcd(a, m)\sim gcd(b, m)$ which is equivalent to $gcd([a]_{\sim}, [m]_{\sim})=[gcd(a, m)]_{\sim}=[gcd(b, m)]_{\sim}=gcd([b]_{\sim}, [m]_{\sim})$ iff   $([a]_{\sim}, [b]_{\sim})\in \tau '_{[m]_{\sim}}$.\hfill\openbox

\begin{corollary}\label{cr100} Let $D$ be a unique factorization domain. Then, for every non zero element $[m]_{\sim}$ of $D'=D_{mult}/\sim$, we have $\tau '_{[m]_{\sim}}=P_{J([m]_{\sim})}$ and so the factor semigroup $D'/\tau '_{[m]_{\sim}}$ satisfies Condition $(*)$ of Definition~\ref{df1}.
\end{corollary}

\noindent
{\bf Proof}. As $\sim \subseteq P_{J(m)}=\tau _m$ by Lemma~\ref{simsub} and Theorem~\ref{th3}, we have $P_{J(m)}/\sim=\tau _m/\sim$ in $D'$.
By Lemma~\ref{ttt}, $P_{J(m)}/\sim=P_{J([m]_{\sim})}$. By Lemma~\ref{sss}, $\tau _m/\sim=\tau '_{[m]_{\sim}}$. Thus $\tau '_{[m]_{\sim}}=P_{J([m]_{\sim})}$ and so the factor semigroup $D'/\tau '_{[m]_{\sim}}$ satisfies Condition $(*)$ of Definition~\ref{df1} by Theorem~\ref{th1}.\hfill\openbox.


\bigskip

If $a$ is an arbitrary nonzero element of a unique factorization domain $D$ then $a$ is either a unit of $D$ or $a$
can be written as a product of (finite many) prime elements (or irreducible elements), uniquely up to order and units.
Let $d(a)$ denote the number of all non associated divisors of $a$, that is, the number of all pairwise different divisors of $[a]_{\sim}$ in $D'=D_{mult}/\sim$. For example, if $D$ is a polynomial ring ${\mathbb Q}[x]$ over the field $\mathbb Q$ of all rational numbers and $a=f(x)=x^3+2x^2-x-2$ then $f(x)=(x-1)(x+1)(x+2)$ and so the pairwise non associated monic divisors of $f(x)$ are
$1$, $x-1$, $x+1$, $x+2$, $(x-1)(x+1)$, $(x-1)(x+2)$, $(x+1)(x+2)$, $(x-1)(x+1)(x+2)$. Thus $d(f(x))=8$. We note that if $a$ is a unit of a unique factorization domain $D$ then $d(a)=1$.

\begin{theorem}\label{eucdx} If $a$ is a nonzero element of a unique factorization domain $D$ then the index of the congruence $P_{J([a]_{\sim})}$ in the semigroup $D'=D_{mult}/\sim$ is finite, and $d(a)=|D'/P_{J([a]_{\sim})}|$.
\end{theorem}

\noindent
{\bf Proof}. Let $a$ be a nonzero element of a unique factorization domain $D$. For an arbitrary divisor $c$ of $a$,
\[gcd([c]_{\sim}, [a]_{\sim})=[gcd(c, a)]_{\sim}=[c]_{\sim}.\] Thus if $x$ and $y$ are non associated divisors of $a$ then, in $D'=D_{mult}/\sim$,
\[([x]_{\sim}, [y]_{\sim})\notin \tau '_{[a]_{\sim}}.\]
Let $c$ be an arbitrary element of $D$. Then \[gcd(gcd(c, a), a)\sim gcd(c, a),\] because $gcd(c, a)$ is a divisor of $a$. Thus
\[(gcd(c, a), c)\in \tau _{a}\] and so \[([gcd(c, a)]_{\sim}, [c]_{\sim})\in \tau _a/\sim=\tau '_{[a]_{\sim}}.\] As $[gcd(c, a)]_{\sim}=gcd([c]_{\sim}, [a]_{\sim})$, we have
$(gcd([c]_{\sim}, [a]_{\sim}), [c]_{\sim})\in \tau '_{[a]_{\sim}}$
from this it follows that $gcd([c]_{\sim}, [a]_{\sim})$ belongs to the $\tau '_{[a]_{\sim}}$-class of $D'$ containing $[c]_{\sim}$. This and the above result together imply that every $\tau '_{[a]_{\sim}}$-class of $D'$ contains one and only one divisor of $[a]_{\sim}$. Consequently, using also Corollary~\ref{cr100}, we have
\[d(a)=|D'/\tau '_{[a]_{\sim}}|=|D'/P_{J([a]_{\sim})}|.\]\hfill\openbox

\section{Applications}

In this section we apply our results for special unique factorization domains.

\begin{theorem}\label{pol} For a nonzero polynomial $f(x)$ of a polynomial ring
${\mathbb F}[x]$ over a field $\mathbb F$, $d(f(x))=|({\mathbb F}[x])'/P_{J([f(x)]_{\sim})}|$, where $({\mathbb F}[x])'={\mathbb F}[x]_{mult}/\sim$.
\end{theorem}

\noindent
{\bf Proof}. By Theorem~\ref{eucdx}, it is obvious.\hfill\openbox

\begin{theorem}\label{positive} For an arbitrary positive integer $m$, $d(m)=|N/P_{J(m)}|$, where $N$ is the multiplicative semigroup of all positive integers and $P_{J(m)}$ is the principal congruence on $N$ defined by the ideal $J(m)$ of $N$ generated by $m$.
\end{theorem}

\noindent
{\bf Proof}. The ring $\mathbb Z$ of all integers is a unique factorization domain. The factor semigroup ${\mathbb Z}'={\mathbb Z}_{mult}/\sim$ is isomorphic to the multiplicative semigroup $N\cup \{ 0\}$ of all non negative integers. Let $m$ be a positive integer. Then, by Theorem~\ref{eucdx},
$d(m)=|({\mathbb Z})'/P_J|$, where $J$ denotes the ideal of $N\cup \{ 0\}$ generated by $m$. By Remark~\ref{rmttt},
$(N\cup \{ 0\})/P_J\cong N/P_J(m)$ and so $d(m)=|N/P_{J(m)}|$.\hfill\openbox

\begin{example}\label{ex2}
Let $S=\{ 0, 1, 2, 3\}$ be the semigroup defined by Table~\ref{Fig2}.

\begin{table}[ht]
 \center{\begin{tabular}{l|l l l l }
 &$1$&$2$&$3$&$0$\\ \hline
$1$&$1$&$2$&$3$&$0$\\
$2$&$2$&$2$&$0$&$0$\\
$3$&$3$&$0$&$3$&$0$\\
$0$&$0$&$0$&$0$&$0$\\
 \end{tabular}}
\caption{}\label{Fig2}
\end{table}

$S$ is a semilattice, that is, a commutative semigroup in which every element is an idempotent.
It is a matter of checking to see that $S\cong N/P_{J(6)}$ and so $d(6)=4$ by Theorem~\ref{positive}.
\end{example}


A quadratic number field is a field of the form $Q[\sqrt d]$, where $d$ is a fixed integer,
positive or negative, which is not a square in $Q$. Its elements are the complex numbers
$a + b\sqrt d$, with $a$ and $b$ in $Q$. The field $Q[\sqrt d]$ is a real quadratic number field if $d > 0$, and
an imaginary quadratic number field if $d < 0$.

A complex number $\alpha$ is called an algebraic number if there is a polynomial $p(x)$ with rational coefficients such that $p(\alpha )=0$.
The monic irreducible polynomial $p(x)$ which generates the kernel of the substitution homomorphism $f(x)\to f(\alpha)$ ($f(x)\in \mathbb Q[x]$) is called the irreducible polynomial for $\alpha$ over $\mathbb Q$.
An algebraic number is called an algebraic integer if its (monic) irreducible polynomial over Q has integer coefficients.
By Proposition 13.1.6, the algebraic integers of a quadratic number field form a ring.

\begin{theorem} Let $d$ be one of the integers $-1$, $-2$, $-3$, $-7$, $-11$, $-19$, $-43$, $-67$, $-163$. If $A$ is a nonzero ideal of the ring $R$ of algebraic integers of an imaginary quadratic number field ${\mathbb Q}[\sqrt d]$ then there is a nonzero element $m$ of $R$ such that $P_A=\tau _m$. Moreover, there is a positive ordinary integer $n$ such that $P_{A\overline A}=\tau _n$, where $\overline A=\{ \overline a:\ a\in A\}$, the complex conjugate of $A$.
\end{theorem}

\noindent
{\bf Proof}.  As $d$ is one of the integers $-1$, $-2$, $-3$, $-7$, $-11$, $-19$, $-43$, $-67$, $-163$, the ring $R$ of algebraic integers of an imaginary quadratic number field is a unique factorization domain by Theorem 13.2.5 of \cite{Artin}. Then, by Theorem 13.5.6 of \cite{Artin}, $R$ is a principal ideal domain. Thus, if $A$ is a nonzero ideal of $R$ then there is a nonzero element $m$ of $R$ such that $I=J(m)=mR$. By Theorem~\ref{th3}, we get $P_A=P_{J(m)}=\tau _m$. By Lemma 13.4.8 of \cite{Artin}, there is a positive ordinary integer $n$ such that $A\overline A=J(n)$ and so $P_{A\overline A}=P_{J(n)}=\tau _n$ by Theorem~\ref{th3} of our paper.\hfill\openbox.

\bigskip

\noindent
Attila Nagy

\noindent
Department of Algebra

\noindent
Mathematical Institute

\noindent
Budapest University of Technology and Economics

\noindent
e-mail: nagyat@math.bme.hu

\end{document}